\documentclass[12pt]{amsart}
\usepackage{mathtools}
\usepackage{amssymb}

\renewcommand{\(}{\left(}
\renewcommand{\)}{\right)}

\newcommand{\<}{\langle}
\renewcommand{\>}{\rangle}
\newcommand{\x}{\times}
\renewcommand{\bar}{\overline}
\newcommand{\abs}[1]{\left\lvert#1\right\rvert}
\newcommand{\norm}[1]{\left\lVert#1\right\rVert}
\newcommand{\st}{\:|\:}

\newcommand{\CC}{{\mathbb{C}}}
\newcommand{\RR}{{\mathbb{R}}}
\newcommand{\ZZ}{{\mathbb{Z}}}

\renewcommand{\phi}{\varphi}
\newcommand{\p}{\rho}

\renewcommand{\span}{{\mathrm{span}}}

\newcommand{\tr}{{\mathrm{tr}}}

\newcommand{\B}{{\mathcal{B}}}

\renewcommand{\H}{{\mathcal{H}}}

\newcommand{\BH}{{\mathcal{B}}(\H)}
\newcommand{\A}{{\mathcal{A}}}

\renewcommand{\tilde}{\widetilde}

\theoremstyle{plain}
\newtheorem{thm}{Theorem}[section]
\newtheorem{lem}[thm]{Lemma}
\newtheorem{prop}[thm]{Proposition}

\theoremstyle{definition}
\newtheorem{defn}[thm]{Definition}

\theoremstyle{remark}

\title{A new proof of {Benedicks'} Theorem for the Weyl Transform}

\author{M.~K.~Vemuri}
\address{Department of Mathematical Sciences, IIT (BHU), Varanasi 221 005}


\email{mkvemuri@gmail.com}

\date{August 16, 2018}

\begin{document}

\begin{abstract}
Benedicks theorem for the Weyl Transform states:
If the set of points where a function is nonzero is of finite measure, and
its Weyl transform is a finite rank operator, then the function is identically
zero.  A new, more transparent proof of this theorem is given.
\end{abstract}

\maketitle

%\tableofcontents\newpage

\section{Introduction}\label{S:intro}
That a nonzero function and its Fourier transform cannot both be sharply
localized is known as the {\em Uncertainty Principle} in Harmonic Analysis.
There are many different precise formulations of this principle, depending
on the way in which localization is quantified.  One such is Benedicks'
theorem \cite{Benedicks}: if $f \in L^1(\RR)$, and the sets
$\{x\in\RR \st f(x) \ne 0\}$ and $\{\xi\in\RR \st \hat{f}(\xi) \ne 0\}$
both have finite Lebesgue measure, then $f \equiv 0$.

Let $\H=L^2(\RR)$, and $\B(\H)$ the set of bounded operators on $\H$.
If $f \in L^1(\RR^2)$, the {\em Weyl transform} of $f$ is the operator
$W(f) \in \B(\H)$ defined by
\begin{equation*}
(W(f)\phi)(t) = \iint f(x,y) e^{\pi i (xy+2yt)} \phi(t+x) \, dx dy.
\end{equation*}

The following analogue of Benedicks' theorem for the Weyl transform was
proved in \cite{btwt}.
\begin{thm}\label{T:main-theorem}
If the set $\{w \in \RR^2 \st f(w) \ne 0\}$ has finite Lebesgue measure and
$W(f)$ is a finite rank operator, then $f \equiv 0$.
\end{thm}

Recall that the Heisenberg group $G$ is the set of triples
\begin{equation*}
\{(x,y,z) \st x,y \in \RR, z \in \CC, \abs{z}=1\}
\end{equation*}
with multiplication defined by
\begin{equation*}
(x,y,z)(x',y',z')=\(x+x',y+y', zz'e^{\pi i(xy'-yx')}\).
\end{equation*}
According to the {\em Stone-von Neumann Theorem}, there is a unique
irreducible unitary representation $\p$ of $G$ such that
\begin{equation*}
\p(0,0,z)= zI.
\end{equation*}
The standard realization of this representation is on the Hilbert space
$\H$ by the action
\begin{equation*}
(\p(x,y,z)\phi)(t)=ze^{\pi i(xy+2yt)}\phi(t+x).
\end{equation*}
Thus, the Weyl transform may be expressed as
\begin{equation*}
W(f) = \iint f(x,y) \p(x,y,1) \, dx dy, \qquad f \in L^1(\RR^2).
\end{equation*}

If $X$ is a trace class operator on $\H$, the modified Fourier-Wigner transform
of $X$ is the function $\alpha(X):\RR^2 \to \CC$ defined by
\begin{equation*}
\alpha(X)(x,y)=\tr(X\p(x,y,1)^*).
\end{equation*}

It is well known (see e.g.\ \cite{Folland}) that if
$f \in L^1(\RR^2)$ and $W(f)$ is a trace class
operator then $\alpha(W(f))=f$, and that if $X$ is a trace class operator
on $\H$ and $\alpha(X) \in L^1(\RR^2)$ then $W(\alpha(X))=X$.  Thus we may
reformulate Theorem \ref{T:main-theorem} as

\begin{thm}\label{T:reformulation}
If $X$ is a finite rank operator on $\H$ and the set
$\{w \in \RR^2 \st \alpha(X)(w) \ne 0\}$ has finite measure, then $X=0$.
\end{thm}

Theorem \ref{T:reformulation} was proved in \cite{btwt}.  In spirit,
the proof in \cite{btwt} was similar to that of Benedicks
\cite{Benedicks}.  However, the definition of the periodization map
was opaque, and embedded in messy representation theory.  In this
work, we define the periodization map in a more transparent way, and
give a proof of Theorem \ref{T:reformulation} that is even closer to that of
Benedicks \cite{Benedicks}.  The price for this is the use of somewhat
more sophisticated ideas from the theory of von Neumann algebras.
However, this is a small price to pay, because ultimately,
Theorem \ref{T:reformulation} depends on Linnel's theorem \cite{Linnel},
which uses, in an essential way, very deep ideas from the theory of
von Neumann algebras.

\section{Periodic Operators}\label{S:periodic-operators}

Let $\pi:G \to \RR^2$ be the projection $\pi(x,y,z)=(x,y)$.  Then
$\pi$ is a homomorphism and $\ker(\pi)=Z(G)$, the center of $G$.
Let $U(1)=\{z \in \CC \st \abs{z}=1\}$.  For $(x,y) \in \RR^2$,
let $s(x,y)=(x,y,1)$.  Then $s$ is a section of $\pi$, i.e.
$\pi\circ s=\mathrm{id}_{\RR^2}$, and
\begin{equation*}
s(x,y) s(x',y')=(x+x',y+y',e^{\pi i (xy'-yx')})=\psi((x,y),(x',y'))s(x+x',y+y'),
\end{equation*}
where $\psi((x,y),(x',y')) \in Z(G)$ is defined by
\begin{equation*}
\psi((x,y),(x',y'))=(0,0,e^{\pi i (xy'-yx')}).
\end{equation*}
Define
$e:\RR^2 \x \RR^2 \to U(1)$ by
\begin{equation*}
e((x,y), (x',y'))=e^{2\pi i (xy'-yx')}.
\end{equation*}
Note that if $g,g' \in G$ then
\begin{equation*}
gg'g^{-1}g'^{-1}=(0,0,e(\pi(g),\pi(g'))).
\end{equation*}
The pairing $e$ is an alternating bicharacter, and it is perfect in the sense
that it establishes the Pontryagin duality $\widehat{\RR^2} \cong \RR^2$.

Let $N \subseteq \RR^2$ be a lattice (necessarily cocompact).
Let $N^\perp =\{w \in \RR^2 \st e(n,w)=1 \quad \forall n \in N\}$.

\begin{defn}
A closed, densely defined operator $T: \H \to \H$ is said to be
$N$-periodic if for all $n \in N$, we have
\begin{equation*}
\p(n,1) T \p(n,1)^{-1}= T.
\end{equation*}
\end{defn}
The set $\A_N=\{T \in \BH \st \text{$T$ is $N$-periodic}\}$ is a
von Neumann algebra (it is the commutant of the self-adjoint set
$\{\p(n,1) \st n \in N\}$). Observe that a closed densely defined
operator $T$ is $N$-periodic iff $T$ is affiliated to $\A_N$.

Observe that for each $n'\in N^\perp$, the operator $\p(n',1)$ is
$N$-periodic, hence the weak closure of
$\span\{\p(n',1) \st n' \in N^\perp\}$ is contained in $\A_N$.  The
following theorem is a direct consequence of \cite[Theorem 5.7]{rcr}.
\begin{thm}\label{T:generating-set}
$\A_N$ is the weak closure of $\span\{\p(n',1) \st n' \in N^\perp\}$.
\end{thm}

In \cite{Linnel}, Linnel constructed a faithful finite weakly continuous
tracial state $\tau$ on $\A_N$ with the property that
\begin{equation}\label{E:trig-fns}
\tau(\p(n',1))=
\begin{cases}
1 & \text{if $n'=0$}\\
0 & \text{if $n' \in N^\perp\setminus\{0\}$,}
\end{cases}
\end{equation}
thus proving the following theorem.

\begin{thm}\label{T:finite}
$\A_N$ is a finite von Neumann algebra.
\end{thm}

Let $L^p(\A_N, \tau)$ denote the Dixmier-Segal ``non-commutative''
$L^p$-space of $\A_N$ with respect to the trace $\tau$ (see
\cite{Pisier-Xu}).  Note that $L^2(\A_N, \tau)$ is a Hilbert space
under the inner product $\<A, B\> = \tau(AB^*)$, and we have the
duality $\A_N=L^\infty(\A_N, \tau) \cong L^1(\A_N, \tau)^*$ (under
the same pairing).  Moreover, since $\A_N$ is finite, we have
$\A_N \subseteq L^2(\A_N, \tau) \subseteq L^1(\A_N, \tau)$.

\begin{defn}
If $T \in L^1(\A_N, \tau)$ and $n' \in N^\perp$, we define the
Fourier-Wigner coefficients of $T$ by
\begin{equation*}
(\alpha(T))(n') = \tau(T\p(n',1)^*) \qquad n' \in N^\perp.
\end{equation*}
\end{defn}

\begin{prop}\label{P:FWC-uniqueness}
If $T \in L^1(\A_N, \tau)$ and $\alpha(T)(n')=0$ for all $n'\in N^\perp$, then
$T=0$.
\end{prop}

\begin{proof}
This follows immediately from Theorem \ref{T:generating-set} and the
duality $\A_N \cong L^1(\A_N, \tau)^*$.
\end{proof}

\begin{prop}\label{P:FWC-parseval}
$\{\p(n',1)\}_{n'\in N^\perp}$ is a complete orthonormal set in
$L^2(\A_N, \tau)$.
\end{prop}

\begin{proof}
This follows immediately from equation (\ref{E:trig-fns}) and
Theorem \ref{P:FWC-uniqueness}.
\end{proof}

\section{The Zak transform}\label{S:zak}
As the Zak (Weil-Brezin) transform is used in several proofs in Section
\ref{S:periodization}, we review some of its properties; they are well
known (see e.g. \cite{Daubechies}).

If $\phi \in \H$ is continuous and compactly supported, we define its
Zak transform
$Z\phi: \RR^2 \to \CC$ by
\begin{equation*}
(Z\phi)(\theta,\sigma) = \sum_{l \in \ZZ} \phi(\sigma-l) e^{2\pi il\theta}.
\end{equation*}
Then,
\begin{equation*}
\begin{aligned}
(Z\phi)(\theta+1,\sigma) =&\; (Z\phi)(\theta,\sigma) \quad\text{and}\\
(Z\phi)(\theta,\sigma+1) =&\; e^{2\pi i \theta} (Z\phi)(\theta,\sigma),
\end{aligned}
\end{equation*}
so $(Z\phi)$ is determined, by its values on $\Omega=[0,1)\x [0,1)$.
Also, if $e_{mn}(t)=\chi_{[0,1)}(t-n) e^{2\pi i mt}$, then
\begin{equation*}
(Ze_{mn})(\theta,\sigma)=e^{-2\pi i n\theta} e^{2\pi i m\sigma},
\qquad (\theta, \sigma) \in \Omega
\end{equation*}
so $Z$ extends to a unitary map $\H \to L^2(\Omega)$.  Also, if
$F:\Omega\to\CC$ is continuous, then
\begin{equation*}
(Z^{-1}F)(t)=\int_0^1 F(\theta,t) \, d\theta.
\end{equation*}

Observe that
\begin{equation*}
\begin{aligned}
(Z\p(x,y,z)\phi)(\theta,\sigma)
=&\; \sum_{l\in\ZZ} (\p(x,y,z)\phi)(\sigma-l) e^{2\pi i l\theta}\\
=&\; \sum_{l\in\ZZ} z e^{\pi i (xy+2y(\sigma-l)} \phi(\sigma-l+x) e^{2\pi i l\theta}\\
=&\; z e^{\pi i (xy+2y\sigma)} \sum_{l\in\ZZ} e^{2\pi il(\theta-y)}\phi(\sigma+x-l)\\
=&\; z e^{\pi i (xy+2y\sigma)} (Z\phi)(\theta-y,\sigma+x).
\end{aligned}
\end{equation*}

Let $M=\ZZ^2 \subseteq \RR^2$.  Then $M$ is a maximal isotropic subgroup.
If $(m_1, m_2)\in M$, then
\begin{equation*}
\begin{aligned}
(Z\p(m_1,m_2,e^{\pi i m_1m_2})\phi)(\theta,\sigma)
=&\; e^{\pi i m_1m_2} e^{\pi i (m_1m_2+2m_2\sigma)} (Z\phi)(\theta-m_2,\sigma+m_1)\\
=&\; e^{2\pi i m_2\sigma} e^{2\pi i m_1\theta} (Z\phi)(\theta,\sigma)
\end{aligned}
\end{equation*}

\section{Periodization}\label{S:periodization}

Although most of the results in this section may be formulated and proved
for arbitrary {\em isotropic} lattices, we will henceforth assume
that $N=\ZZ\times a\ZZ$, for some $a\in\ZZ$.  The reason for this is that
it suffices for our main goal, namely the proof of Theorem
\ref{T:reformulation} and that Linnel's trace $\tau$ has a nice formula
in this case.

Observe that $N^\perp=\frac{1}{a}\ZZ\times \ZZ$.  Let
$\chi_h = \chi_{[h/a,(h+1)/a)}$.  If $T \in \A_N$,
then (\cite[p3274]{Linnel})
\begin{equation*}
\tau(T)=\sum_{h=0}^{a-1} \<T\chi_h, \chi_h\>
       =\sum_{h=0}^{a-1} \int_{h/a}^{(h+1)/a} T\chi_h.
\end{equation*}

\begin{lem}
Let $g\in L^1(\Omega)$.  Let
$\mathcal{D}_g=\{\phi\in\H\st g Z\phi \in L^2(\Omega)\}$.
Define $M_g:\mathcal{D}_g \to \H$ by
\begin{equation*}
M_g(\phi) = Z^{-1}(gZ\phi).
\end{equation*}
Then $\mathcal{D}_g$ is dense in $\H$, and $M_g$ is closed.  Moreover
$M_g \in L^1(\A_N, \tau)$.
\end{lem}

\begin{proof}
An elementary measure theory argument shows that $\mathcal{D}_g$ is dense
in $\H$.  Observe that $\mathcal{D}_{\bar{g}}=\mathcal{D}_g$.  Suppose
$\phi_n \in \mathcal{D}_g$, $\phi_n \to \phi \in \H$, and
$M_g\phi_n \to \psi \in \H$.  Then for all $\mu\in\mathcal{D}_{\bar{g}}$, we have
\begin{equation*}
\begin{aligned}
\<\psi,\mu\>
=&\; \lim_{n\to\infty} \<M_g\phi_n, \mu\>\\
=&\; \lim_{n\to\infty} \<\phi_n, M_{\bar{g}}\mu\>\\
=&\; \<\phi, M_{\bar{g}}\mu\>\\
=&\; \<M_g\phi, \mu\>.
\end{aligned}
\end{equation*}
Since $\mathcal{D}_{\bar{g}}$ is dense in $\H$, it follows that $M_g\phi=\psi$.
This proves that $M_g$ is closed.

Extend $g$ to $\RR^2$ by $N$-periodicity.
Then $gZ\phi$ has the same periodicity properties as $Z\phi$.

We claim that $M_g$ is $N$-periodic. Indeed, if $n=(n_1,n_2)\in N$, and
$\phi \in \mathcal{D}_g$, then $Z\p(n,1)\phi \in \mathcal{D}_g$, and
\begin{equation*}
\begin{aligned}
(\p(n,1)^{-1}M_g\p(n,1)\phi)(t)
=&\; (\p(-n,1)Z^{-1}(gZ\p(n,1)\phi))(t)\\
=&\; e^{\pi i (n_1n_2 - 2n_2 t)} (Z^{-1}(gZ\p(n,1)\phi))(t-n_1)\\
=&\; e^{\pi i (n_1n_2 - 2n_2 t)} \int_0^1 (gZ\p(n,1)\phi)(\theta, t-n_1) \, d\theta
\\
=&\; e^{\pi i (n_1n_2 - 2n_2 t)}
     \int_0^1 g(\theta, t-n_1)(Z\p(n,1)\phi)(\theta,t-n_1) \, d\theta\\
=&\; e^{\pi i (n_1n_2 - 2n_2 t)} \int_0^1 g(\theta, t-n_1)
     e^{\pi i (n_1n_2 + 2n_2(t-n_1))} (Z\phi)(\theta-n_2, t)\, d\theta\\
=&\; \int_0^1 g(\theta, t-n_1) (Z\phi)(\theta-n_2, t)\, d\theta\\
=&\; \int_0^1 g(\theta, t) (Z\phi)(\theta, t)\, d\theta\\
=&\; Z^{-1}(gZ\phi)\\
=&\; M_g(\phi).
\end{aligned}
\end{equation*}

Since $\mathcal{A}_N$ is finite, $M_g$ is automatically measurable.  Let
$\lambda>0$.  Observe that $\abs{M_g}=M_{\abs{g}}$, and
$\chi_{(\lambda,\infty)}(\abs{M_g})=M_{g_\lambda}$, where
\begin{equation*}
g_\lambda(\theta,\sigma)=
\begin{cases}
1, & \text{if $\abs{g(\theta,\sigma)}>\lambda$},\\
0, & \text{otherwise.}
\end{cases}
\end{equation*}
Therefore
\begin{equation*}
\begin{aligned}
\tau(\chi_{(\lambda,\infty)}(\abs{M_g}))
=&\; \sum_{h=0}^{a-1} \int_{h/a}^{(h+1)/a} \chi_{(\lambda,\infty)}(\abs{M_g}) \chi_h\\
=&\; \sum_{h=0}^{a-1} \int_{h/a}^{(h+1)/a} Z^{-1}(g_\lambda Z\chi_h)\\
=&\; \sum_{h=0}^{a-1} \int_0^1 g_\lambda(\theta,\sigma) \chi_h(\sigma)
                   \,d\theta \,d\sigma\\
=&\; \int_\Omega g_\lambda(\theta,\sigma) \,d\theta \,d\sigma\\
=&\; m(\{(\theta, \sigma) \st \abs{g(\theta,\sigma)} > \lambda\}),
\end{aligned}
\end{equation*}
i.e., the distribution function of $M_g$ is equal to the distribution
function of $g$.  Since $g\in L^1(\Omega)$, it follows that
$M_g \in L^1(\mathcal{A}_N, \tau)$.

\end{proof}

Let $X=\phi\otimes\bar{\psi} \in \BH$ be a rank-one operator.
%{\color{orange}
Define
\begin{equation*}
\tilde{X}=\sum_{j=0}^{a-1} M_{g_j} \p(-j/a,0,1),
\end{equation*}
where
\begin{equation*}
g_j(\theta,\sigma) = (Z\phi)(\theta,\sigma) \bar{(Z\psi)(\theta,\sigma-j/a)}.
\end{equation*}
%}
Then $\tilde{X} \in L^1(\A_N, \tau)$ by the non-commutative H\"older inequality.
Moreover, the map $(\phi,\psi) \mapsto \widetilde{\phi\otimes\bar{\psi}}$ is
sesquilinear.  If $X$ is a finite rank operator, then there exist
rank-one operators $X_1, \dots, X_n$ such that $X=X_1+\cdots+X_n$.
Define $\tilde{X}=\tilde{X_1}+\cdots+\tilde{X_n}$.  This is independent
of the expression of $X$ in terms of rank-one operators, by the aforementioned
sesquilinearity, and $\tilde{X} \in L^1(\A_N, \tau)$.

%{\color{orange}
\begin{lem}
Let $X\in\BH$ be of the form $X=\phi\otimes\bar{\psi}$ where
$\phi,\psi$ are smooth compactly supported functions.  If $f$ is a
smooth compactly supported function, then
\begin{equation*}
\sum_{n\in N} \p(n,1) X \p(n,1)^{-1} f = a^{-1}\tilde{X}f,
\end{equation*}
in the sense that the square partial sums of the left hand side converge
in $L^2(\RR)$ to the right hand side.
\end{lem}

\begin{proof}
For $L\in\ZZ$, let
\begin{equation*}
X_L=\sum_{k=-L}^L\sum_{l=-L}^L \p(k,la,1)X\p(k,la,1)^{-1}
\end{equation*}
Then 
\begin{equation*}
\begin{aligned}
(ZX_Lf)(\theta,\sigma)
=&\; (ZX_LZ^{-1}Zf)(\theta,\sigma)\\
=&\; \sum_{k=-L}^L\sum_{l=-L}^L
\(\(Z\p(k,la,1)\phi\otimes\bar{Z\p(k,la,1)\psi}\)(Zf)\)(\theta,\sigma)\\
=&\; \sum_{k=-L}^L\sum_{l=-L}^L
\<Zf,Z\p(k,la,1)\psi\> \(Z\p(k,la,1)\phi\)(\theta,\sigma)\\
=&\; \sum_{k=-L}^L\sum_{l=-L}^L e^{2\pi i (la(\sigma-\sigma')+k(\theta-\theta'))}
\iint_\Omega (Zf)(\theta',\sigma')\bar{Z\psi(\theta',\sigma')}\,d\theta'd\sigma'
(Z\phi)(\theta,\sigma)\\
=&\; \iint_\Omega \(\sum_{k=-L}^L\sum_{l=-L}^L e^{2\pi i (la(\sigma-\sigma')+k(\theta-\theta'))}\)
(Zf)(\theta',\sigma')\bar{Z\psi(\theta',\sigma')}
\,d\theta'd\sigma'  (Z\phi)(\theta,\sigma)\\
=&\; \iint_\Omega D_L(a(\sigma-\sigma'))D_L(\theta-\theta')
(Zf)(\theta',\sigma')\bar{Z\psi(\theta',\sigma')}
\,d\theta'd\sigma' (Z\phi)(\theta,\sigma)\\
=&\; \int_0^1 \int_0^1 D_L(a\sigma')D_L(\theta')
(Zf)(\theta-\theta',\sigma-\sigma')\bar{Z\psi(\theta-\theta',\sigma-\sigma')}
\,d\theta'd\sigma' (Z\phi)(\theta,\sigma)\\
=&\; \frac{1}{a} \int_0^a \int_0^1 D_L(\sigma')D_L(\theta')
(Zf)(\theta-\theta',\sigma-\sigma'/a)
\bar{Z\psi(\theta-\theta',\sigma-\sigma'/a)}
\,d\theta'd\sigma' (Z\phi)(\theta,\sigma)\\
=&\; \frac{1}{a} \sum_{j=0}^{a-1} \int_0^1 \int_0^1 D_L(\sigma') D_L(\theta')
(Zf)(\theta-\theta',\sigma-(\sigma'+j)/a)
\bar{Z\psi(\theta-\theta',\sigma-(\sigma'+j)/a)}
\,d\theta'd\sigma'\\
 &\; \x (Z\phi)(\theta,\sigma)\\
\end{aligned}
\end{equation*}
Since $Zf$ and $Z\psi$ are smooth, it follows that
\begin{equation*}
(ZX_Lf)(\theta,\sigma) \to \frac{1}{a} \sum_{j=0}^{a-1}
(Zf)(\theta, \sigma-j/a)\bar{Z\psi(\theta, \sigma-j/a)}(Z\phi)(\theta,\sigma)
=\frac{1}{a}Z\tilde{X}f
\end{equation*}
uniformly, and so in $L^2(\Omega)$.
\end{proof}

\begin{lem}\label{L:local-rank}
If $\operatorname{rank}(X) < a$ then $\tilde{X}$ is not injective.
\end{lem}

We will need the following measure theoretic lemma to prove Lemma
\ref{L:local-rank}.

\begin{lem}\label{L:measure-theoretic-lemma}
Let $\Sigma$ be a measurable space, $m < n$, and
$e_1, \dots, e_m:\Sigma \to \CC^n$ be measurable.  Then there exists
$e: \Sigma \to \CC^n$ measurable such that
$\norm{e(x)}=1$, and $\<e_j(x), e(x)\>=0$ for all $x \in X$, $j=1,\dots,m$.
\end{lem}

\begin{proof}[Proof of Lemma \ref{L:local-rank}]
Let $b=\operatorname{rank}(X)$, and assume $b<a$.  Then there exist
$\phi_j, \psi_j \in \H$, $j=1,\dots, b$, such that
$X=\phi_1\otimes\bar{\psi_1}+\cdots+\phi_b\otimes\bar{\psi_b}$.
Let $\Sigma=[0,1)\x [0,1/a)$, and for $j=1,\dots, b$,
define $e_j:\Sigma\to\CC^a$ by
\begin{equation*}
e_j(\theta,\sigma)=\(
                   Z\psi_j(\theta,\sigma),
                   Z\psi_j(\theta,\sigma+1/a),
                   \dots,
                   Z\psi_j(\theta,\sigma+(a-1)/a)
                   \).
\end{equation*}
Therefore, by Lemma \ref{L:measure-theoretic-lemma}, there exists
$e:\Sigma\to\CC^a$ such that $\norm{e(\theta,\sigma)}=1$, and
\begin{equation*}
\sum_{i=1}^a \bar{Z\psi_j(\theta,\sigma+(i-1)/a)}(e(\theta,\sigma))_i = 0,
\qquad (\theta,\sigma)\in\Sigma,\, j=1,\dots,b,
\end{equation*}
where
$(\cdot)_i$ denotes the $i^{\mathrm{th}}$ component.

Define $F:\Omega\to\CC$ by
$F(\theta,\sigma)=(e(\theta, \sigma-(i-1)/a))_i$ if
$\sigma\in[(i-1)/a,i/a)$, $i=1,\dots,a$, where
$(\cdot)_i$ denotes the $i^{\mathrm{th}}$ component.
Clearly $F\in L^2(\Omega)$.
Put $f=Z^{-1}F$.  Then $f\ne 0$. Let
$g_{ji}(\theta,\sigma)=(Z\phi_j)(\theta,\sigma)\bar{Z\psi_j(\theta,\sigma+i/a)}$.
Then
\begin{equation*}
\begin{aligned}
((Z\tilde{X})f)(\theta,\sigma)
=&\; \sum_{j=1}^b \(Z\(\widetilde{\phi_j\otimes\bar{\psi_j}}\)f\)(\theta,\sigma)\\
=&\; \sum_{j=1}^b Z\(\sum_{i=0}^{a-1} M_{g_{ji}} \p(i/a,0,1)f\)(\theta,\sigma)\\
=&\; \sum_{j=1}^b \sum_{i=0}^{a-1} (Z\phi_j)(\theta,\sigma)
                               \bar{Z\psi_j(\theta,\sigma+i/a)}
                               F(\theta,\sigma+i/a)\\
=&\; \sum_{j=1}^b (Z\phi_j)(\theta,\sigma) \sum_{i=0}^{a-1}
                               \bar{Z\psi_j(\theta,\sigma+i/a)}
                               F(\theta,\sigma+i/a)\\
=&\; 0, \qquad (\theta,\sigma)\in \Sigma.
\end{aligned}
\end{equation*}
Since $\sum_{i=0}^{a-1}\bar{Z\psi_j(\theta,\sigma+i/a)}F(\theta,\sigma+i/a)$
is $(1/a)$-periodic in $\sigma$, it follows that $Z\tilde{X}f \equiv 0$, and
hence $\tilde{X}f=0$.
\end{proof}
%}

\begin{lem}
Let $k,l \in \ZZ$.  Then
\begin{equation*}
\alpha(M_g)(k/a,l)=
\begin{cases}
e^{\pi i lk/a} \int_\Omega g(\theta,\sigma)
e^{-2\pi i(l\sigma+k\theta/a)} \, d\theta d\sigma,
& \text{if $k \in a\ZZ$},\\
0,
& \text{otherwise.}
\end{cases}
\end{equation*}
\end{lem}

\begin{proof}

Since $L^\infty(\Omega) \cap L^1(\Omega)$ is dense in $L^1(\Omega)$, it
suffices to check this for $g \in L^\infty(\Omega) \cap L^1(\Omega)$.
Then

\begin{equation*}
\begin{aligned}
\alpha(M_g)(k/a,l)
=&\; \tau(M_g\p(-k/a, -l, 1))\\
=&\; \sum_{h=0}^{a-1} \<M_g\p(-k/a,-l,1)\chi_h, \chi_h\>\\
=&\; \sum_{h=0}^{a-1} \<Z^{-1}(gZ\p(-k/a,-l,1)\chi_h), \chi_h\>\\
=&\; \sum_{h=0}^{a-1} \<gZ\p(-k/a,-l,1)\chi_h, Z\chi_h\>\\
=&\; \sum_{h=0}^{a-1} \int_\Omega g(\theta,\sigma)
(Z\p(-k/a,-l,1)\chi_h)(\theta,\sigma) \bar{(Z\chi_h)(\theta,\sigma)}
\, d\sigma d\theta\\
=&\;
\begin{cases}
\sum_{h=0}^{a-1} \int_\Omega g(\theta,\sigma)
e^{\pi i (lk/a-2l\sigma)} (Z\chi_h)(\theta+l,\sigma-k/a)
\bar{(Z\chi_h)(\theta,\sigma)},
& \text{if $k \in a\ZZ$},\\
0,
& \text{otherwise.}
\end{cases}\\
=&\;
\begin{cases}
\sum_{h=0}^{a-1} \int_\Omega g(\theta,\sigma)
e^{\pi i (lk/a-2l\sigma-2k\theta/a)} (Z\chi_h)(\theta,\sigma)
\bar{(Z\chi_h)(\theta,\sigma)},
& \text{if $k \in a\ZZ$},\\
0,
& \text{otherwise.}
\end{cases}\\
=&\;
\begin{cases}
\int_\Omega g(\theta,\sigma)
e^{\pi i (lk/a-2l\sigma-2k\theta/a)},
& \text{if $k \in a\ZZ$}\\
0,
& \text{otherwise.}
\end{cases}
\end{aligned}
\end{equation*}

\end{proof}

%{\color{red}

\begin{thm}\label{T:poisson}
If $X\in\BH$ is a finite rank operator, then
\begin{equation*}
\alpha(\tilde{X})(n') = \alpha(X)(n')
\end{equation*}
for all $n' \in N^\perp$.
\end{thm}

\begin{proof}
By linearity, it suffices to check this assuming $X$ has rank one.
Let $n' \in N^\perp$.  Then $n'=(k/a,l)$ for some $k,l \in \ZZ$.
There exists a unique $k' \in \{0,1,\dots,a-1\}$ such that
$k+k' \in a\ZZ$.
Write $X=\phi\otimes\bar{\psi}$, $\phi,\psi\in\H$.  Then
\begin{equation*}
\begin{aligned}
\alpha(\tilde{X})(n')
=&\; \tau\(\sum_{j=0}^{a-1} M_{g_j} \p(-j/a,0,1) \p(-k/a,-l,1)\)\\
=&\; \sum_{j=0}^{a-1} \tau(M_{g_j} \p(-(j+k)/a, -l, e^{\pi ijl/a}))\\
=&\; \sum_{j=0}^{a-1} e^{\pi ijl/a} \alpha(M_{g_j})((j+k)/a, l)\\
=&\; e^{\pi i k'l/a} \alpha(M_{g_{k'}}((k'+k)/a, l)\\
=&\; e^{\pi ik'l/a} e^{\pi il(k'+k)/a} \int_\Omega g_{k'}(\theta,\sigma)
     e^{-2\pi i (l\sigma +(k'+k)\theta/a)} \, d\theta d\sigma\\
=&\; e^{\pi il(2k'+k)/a}
     \int_\Omega (Z\phi)(\theta,\sigma) \bar{(Z\psi)(\theta, \sigma-k'/a)}
     e^{-2\pi i (l\sigma +(k'+k)\theta/a)} \, d\theta d\sigma\\
=&\; e^{\pi il(k'+k)/a}\int_\Omega (Z\phi)(\theta,\sigma)
     \bar{(Z\psi)(\theta-l,\sigma-k'/a) e^{\pi i (-lk'/a + 2l\sigma + 2(k'+k)\theta/a)}}
     \, d\theta d\sigma\\
=&\; e^{\pi il(k'+k)/a}\int_\Omega (Z\phi)(\theta,\sigma)
     \bar{(Z\p(-k'/a,l,1)\psi)(\theta,\sigma) e^{2\pi i (k'+k)\theta/a}}
     \, d\theta d\sigma\\
=&\; e^{\pi il(k'+k)/a}\int_\Omega (Z\phi)(\theta,\sigma)
     \bar{(Z\p(-k'/a,l,1)\psi)(\theta,\sigma+(k'+k)/a)}
     \, d\theta d\sigma\\
=&\; \int_\Omega (Z\phi)(\theta,\sigma)
     \bar{(Z\p((k'+k)/a,0,e^{-\pi i l(k'+k)/a})\p(-k'/a,l,1)\psi)(\theta,\sigma)}
     \, d\theta d\sigma\\
=&\; \int_\Omega (Z\phi)(\theta,\sigma)
     \bar{(Z\p(k/a,l,1)\psi)(\theta,\sigma)} \, d\theta d\sigma\\
=&\; \<\phi, \p(k/a,l,1)\psi\>\\
=&\; \tr((\phi\otimes\bar{\psi})\p(k/a,l,1)^*)\\
=&\; \alpha(X)(n')
\end{aligned}
\end{equation*}
\end{proof}

%}

\section{The Benedicks argument}\label{S:benedicks}

The proof of Theorem \ref{T:reformulation} proceeds as in \cite{Benedicks,btwt}.
Assume $X$ is a finite rank operator on $\H$.  Let
$B=\{w \in \RR^2 \st \alpha(X)(w) \ne 0\}$ and assume that $B$ has finite
measure.  Choose an integer $a$ greater than the rank of $X$ and let
$N=\ZZ\x a\ZZ$.  Then $\tilde{X}$ is not injective by Lemma \ref{L:local-rank}.
Also $N_v=(B-v) \cap N^\perp$ is a finite set for almost
every $v \in \RR^2$.  For such $v$, let $X^v=X\p(s(v))^*$. Since
$\alpha(X^v)(w) = \bar{\psi(w,v)}\alpha(X)(w+v)$, it follows that
$\alpha(X^v)(w)=0$ if $w \notin B-v$.  Therefore by Theorem \ref{T:poisson},
$\alpha(\tilde{X^v})(n')=\alpha(X^v)(n')=0$ if $n' \in N^\perp \setminus N_v$.
Therefore, by Proposition \ref{P:FWC-uniqueness} and
Proposition \ref{P:FWC-parseval},
\begin{equation*}
\tilde{X^v} = \sum_{n'\in N_v} \alpha(X^v)(n) \p(n,1).
\end{equation*}
Since the right hand side is a finite sum, and $\tilde{X^v}$ is not injective,
it follows that $\tilde{X^v} = 0$ by \cite[Theorem 1.2]{Linnel}.
Therefore $\alpha(X)(n'+v)=0$ for all $n' \in N^\perp$ and almost every
$v \in \RR^2$.  Therefore $\alpha(X)=0$ almost everywhere, and so $X=0$.

\bibliographystyle{amsplain}
\bibliography{v8-npbtwt}

\end{document}